# NECESSARY AND SUFFICIENT CONDITIONS FOR DIFFERENTIABILITY OF A FUNCTION OF SEVERAL VARIABLES.

R.P. VENKATARAMAN,
#1371, 13'TH MAIN ROAD,
II STAGE, FIRST PHASE,
B.T.M. LAYOUT,
BANGALORE 560 076.
<rpvraman@yahoo.com>

## ABSTRACT
The necessary and sufficient conditions for differentiability of a function of several real variables stated and proved and its ramifications discussed.

**Proposition:** A scalar or vector valued function $f$ of $N$ real variables $x_1, x_2, \ldots x_N$ is differentiable at $(c_1, c_2, \ldots c_N)$ iff with respect to the point under consideration as origin in spherical polar coordinates $(r\phi_1, \phi_2, \phi_3, \ldots \phi_{(N-2)}, \vartheta)$ it is differentiable at the origin i.e.

$$df = \frac{\partial f}{\partial r} dr$$

exists, at $r=0$ i.e. it is independent of the angle variables.
Proof.:
The limit of a function $f(x_1, x_2, \ldots x_N)$ at $(c_1, c_2, \ldots c_N)$ is independent of the choice of the basis. Hence in spherical polar coordinates with $(c_1, c_2, \ldots c_N)$ as origin if $lt_{r \to 0} f$ is independent of the angle variables $(\phi_1, \phi_2, \phi_3, \ldots \phi_{(N-2)}, \vartheta)$ then the limit exists and the converse holds by the very definition of limit.

Derivative being a limit the result stated in the proposition follows. The differential of a function of several variables is given by

$$df = \nabla f . dr = \frac{\partial f}{\partial r} dr$$

At $r=0$ it is easy to show that the partial derivatives of $f$ with respect to the angle variables are all zero since the arguments are all zero and hence the difference vanishes..
Hence the proposition.

In the case of vector valued function the particular choice of the basis reduces the division to that by a scalar.

Cor.:1 A function $f$ of complex variable $z$ given by $f(z) = u + iv$ is differentiable at $z_0$ iff in polar coordinates defined by
$$z - z_0 = r \exp(i\vartheta)$$
the derivatives $P_r$ and $Q_r$ of the real and imaginary parts of $f$ exist at $r = 0$.

Cor.: 2 A function $f$ of quaternion variable **h**=a+**i**b+**j**c+**k**d is differentiable at $h_0$ iff in spherical polar coordinates with $h_0$ as origin it is differentiable at the origin.

**Discussion**: The choice of spherical polar coordinates covers all possible paths and helps us check easily if the limit exists and again find the limit easily whenever it exists. In the case of differentiability of $f(z)$ the proposition encompasses Cauchy-Riemann conditions and leads to important conclusions. To check the continuity of the derivative one has to choose a generic point as the origin. The problem of finding the maxima and minima of functions of several variables is drastically simplified to the problem of solving one equation in radial coordinate and checking the signature of the second derivative in that coordinate. These are illustrated below with examples.

It has also equipped one to proceed with calculus on the division ring of quaternions where multiplication is noncommutative.

**Examples:**
Limit at the origin:

1. $xy/[x^2+y^2]$

In the standard method it is shown that along $y=mx$, the limit is $\dfrac{m}{1+m^2}$ and hence it does not exist. In polar coordinates it is $lt_{r \to 0} \sin\Theta\cos\Theta$ and hence it does not exist as the angle variables are undefined at the origin.

2. $y\exp(-1/x^2)/[y^2+\exp(-2/x^2)]$

The limit is 0 along $y=ax^n$ for any $n$ but is ½ along $y=\exp(-1/x^2)$. In polar coordinates it is

$$lt_{r \to 0} \exp(-1/(r^2\cos^2\Theta)/(r\sin\Theta)$$

The above limit is infinity when $\Theta=0$ or $\Pi$ and is 0 otherwise and hence undefined.

3. $x^2y^2/[x^2+y^2]$

In polar coordinates the limit is $lt_{r \to 0}(r\sin\Theta\cos\Theta)^2$ and is found to be 0.

Example 4. Differential of $\sqrt{xy}$

$$df = 1/2\sqrt{y/x}dx + 1/2\sqrt{x/y}dy$$

does not exist at the origin as the partial derivatives do not exist at the origin.
In polar coordinates

$$dF = \sqrt{\sin(2\Theta)/2}dr + r/\sqrt{2}\cos(2\Theta)/\sqrt{\sin(2\Theta)}d\Theta$$

does not exist as the angle is undefined at the origin.

Example 5. Maxima and minima of $f(x,y)$

The procedure is to find the derivative of the function at a generic point *P(a,b)* by shifting the origin to this point using spherical polar coordinates and solve for the ordered pair *(a,b)* by equating the derivative to zero

$$f(x,y) = xy\exp(-xy)$$

$$df = [\{y-xy^2\}\exp(-xy)]dx + [\{x-x^2y\}\exp(-xy)]dy$$

In Cartesian coordinates one finds the solutions to be $x=y=0 \text{ and } xy=1$. From the second derivatives it is easy to show that the second represents a maximum and the

first neither a maximum or a minimum as it does not satisfy $f_{xx}f_{yy} > (f_{xy})^2$. In polar coordinates

$$F = r^2(\sin 2\Theta)/2 \exp[-(r^2/2)\sin(2\Theta)]$$
$$\partial f / \partial r = [r\sin(2\Theta) - r^3(\sin^2 2\Theta)/2]\exp(-r^2 \sin(2\Theta)/2)$$
$\partial f / \partial r = 0$ yields the same solutions as above:
$$r = 0 \text{ and } r^2\{\sin(2\Theta)/2\} = 1.$$

The second derivative $\partial^2 f / \partial r^2$ does not exist at the origin and satisfies

$$\partial^2 f / \partial r^2 < 0 \text{ if } r^2 \sin\Theta\cos\Theta = 1$$

Thus the problem of finding the maxima and minima is reduced to solving just one equation in *r* coordinate and checking the signature of the second derivative in that coordinate.

Example 6. In the examples below a generic point is chosen as the origin and *f(z)* is single

I) Derivative of $f(z) = \bar{z}$

Using $z_0$ as the origin in polar coordinates, we find
$$f'(z) = \exp(-i\theta)[\exp(-i\theta)]$$
and hence it does not exist anywhere.

II) Derivative of $f(z) = |z|^2$

$$f(z) = |z|^2 = (x-x_0)^2 + (y-y_0)^2 + (x_0^2 + y_0^2) + 2(x-x_0)x_0 + 2(y-y_0)y_0$$
Hence $f'(z) = \exp(-i\vartheta)[2r + 2x_0\cos\vartheta + 2y_0\sin\vartheta]$
exists at *(0,0)* only.

III) Derivative of $f(z) = z^n$ where *n* is an integer.
$$f'(z) = \exp(-i\Theta)[nr^{(n-1)}\cos n\Theta + inr^{(n-1)}\sin n\Theta]$$
Since $P_r(r,\Theta)$ and $Q_r(r,\Theta)$ do not exist at *r=0*, if *n< 1*, *f(z)* is not differentiable at *r = 0* whenever *n <1*. It is easy to show as in the above two examples, by taking a generic point as the origin that the derivative in this case is a continuous function of $z_0$ whenever *n≥1*.

IV) Proof of existence of higher derivatives for an analytic function from first principles:
Setting *t=rexp(iθ)* we get

$$f'(z)|_{(z=z_0)} = \exp(-i\vartheta)\,[lt_{(t \to 0)} \frac{d}{dt} f(t+z_0)]\exp(i\vartheta)$$

The second derivative by definition is

$$f''(z)|_{(z=z_0)} = \exp(-i\vartheta)\,[lt_{(\Delta t \to 0)} \frac{f'(\Delta t, z_0) - f'(0, z_0)}{\Delta t}]\exp(i\vartheta)$$

and since analyticity at *z=$z_0$* implies the continuity of *f'(z)*

$$\left|[f''(z)]_{(z=z_0)}\right| < \exp(-i\Theta) \, [lt_{(\delta \to 0)} \frac{\varepsilon(\delta)}{\delta}] \exp(i\Theta)$$

where $\varepsilon$ and $\delta$ are positive quantities. The above inequality shows $f''(z_0)$ is finite. The limit is unique since $r$, and hence $t$, tend to zero only from the positive side. It is easy to show that the second derivative is continuous since $f''(z_0+\Delta z_0)$ tends to $f''(z_0)$ as $\Delta z_0$ tends to zero. Hence, one can similarly show that the higher derivatives also exist. This proof holds even for a function of a real variable, which is defined by one expression in its entire domain, except that the values of $\theta$ are only zero and $\pi$. The same definition of analyticity resulting in the existence of all higher derivatives holds also for a function of several real variables and will be discussed in a subsequent article.

Thus continuity of $f'(x)$ implies the existence of higher derivatives. The proof of existence of Laurent's series based on Fourier series given by Myskis[1], and hence that of Taylor series also hold for $f(x)$. It is only while determining the radius of convergence of the series caution has to be exercised to take into account the singularities in the complex plane. The coefficients are given by the corresponding values for $f(z)$. But in practice it can be found directly using Taylor expansion of known functions. For example $(1+x^2)^k$ where $k$ is real, has a power series with radius of convergence $1$ limited by the singularities on the imaginary axis. It is easy to show that $f(x)$ is analytic at $x=\pm 1$ and has a Taylor series expansion about $x^2=1$ whose radius of convergence is limited by the singularities on the imaginary axis. That $exp(-1/x^2)$ does not have a power series is because $1/x^2$ is not analytic at the origin although exponential function is an entire function. Thus analyticity for a function of real variable, which is defined by one expression in its domain, could be defined in the same way as for functions of complex variable.

The choice of "Spherical Polar Coordinates" for this note would be equally appropriate the result hinging as it does on this specific choice. Even in power series expansion of functions of several variables this choice comes handy in addressing convergence questions and also in the solution of first order differential equations :
The equation

$$dy/dx = f(x,y) = -P(x,y)/Q(x,y)$$

can be solved easily in $(r, \Theta)$ using the method of separation of variables as long as $P$ and $Q$ are homogeneous in the two variables and of the same degree for

$$dr/d\vartheta = \frac{P r \sin\vartheta - Q r \cos\vartheta}{Q \sin\vartheta + P \cos\vartheta}$$

From the above it is easy to realise that singular points be rather defined as those through which more than one curve passes violating Cauchy's theorem. For, in this basis, for the equation governing a family of circles origin is not a singular point, the numerator being zero and the denominator unity.

Taylor expansion in two real variables can be written as[2]

$$f(x,y) = \sum_{(i+j)=0}^{n} a_{ij} x^i y^j + R_n$$

where $R_n$ is the remainder term and if $a_{ij}$ is defined as[3]

$$a_{ij} = \frac{1}{i!\,j!} [\partial^{(i+j)}/\partial x^i \partial y^j f]_{(x=0, y=0)}$$

then Taylor series is defined as the limit of the above sum since convergence could be addressed more easily in polar coordinates. For example

1) $(1+xy)^k = 1 + k(xy)/1! + k(k-1)(xy)^2/2! + \ldots$

is convergent for all *(x,y)εR²* and *kεR* whenever in polar coordinates *r* is less than one.

2) $\exp(xy) = 1 + (xy)/1! + (xy)^2/2! + \ldots + (xy)^r/r! + \ldots$

is convergent for all *(x,y)εR²*. Power series for a few cases are given by Kaplan[4]. Extension to several variables is straight forward.

**Conclusion :**

Using spherical coordinates, necessary and sufficient conditions for the existence of the limit of a scalar / vector valued function of several variables at a point and hence those for differentiability of the above functions at a point and hence those for differentiability of a function of complex variable and also quaternion variable have been stated and proved. The problem of finding the maxima/minima is drastically simplified to solving just one equation in radial coordinate and checking the nature of the second derivative in that coordinate. The conditions for differentiability of *f(z)* encompasses Cauchy-Riemann conditions. It has been proved from first principles that analyticity of a function at a point implies the existence of all higher derivatives and that the same definition is valid even for functions of real variable.. A few examples have been discussed to illustrate the constructive method of finding the limit using the above choice of the basis.

No counter examples could be produced without violating the topological nature of limits. For further study functions of several complex variables, calculus on quaternions, solutions to differential equations about irregular singular points and ramifications ,if any, on variational methods are being taken up.

*************************************************************************